\documentclass[12pt]{article}

\usepackage[monochrome]{color}
\usepackage{a4wide}
\usepackage{amsmath}
\usepackage{amsfonts}
\usepackage{amssymb}

\usepackage{epsfig}
\usepackage{graphicx}

\newcommand{\IR}{\ensuremath{\mathbb{R}}}

\newcommand{\IQ}{\ensuremath{\mathbb{Q}}}

\newcommand{\TT}{\ensuremath{{\cal{T}}}}

\newcommand{\HH}{\ensuremath{\mathcal{H}}}

\newtheorem{theorem}{Theorem}[section]
\newtheorem{lemma}[theorem]{Lemma}
\newtheorem{hypothesis}[theorem]{Hypothesis}
\newtheorem{corollary}[theorem]{Corollary}

\date{}

\begin{document}

\title{
{\bf Trivial, critical and near-critical scaling limits of two-dimensional percolation}}

\author{
{\bf Federico Camia},
\thanks{Research supported in part by a Veni grant of the NWO
(Dutch Organization for Scientific Research). E-mail: fede@few.vu.nl}
{\bf \, Matthijs Joosten},
\thanks{E-mail: mjoosten@few.vu.nl}
{\bf \, Ronald Meester}
\thanks{Research supported in part by a Vici grant of the NWO
(Dutch Organization for Scientific Research). E-mail: rmeester@few.vu.nl}\\
{\sl Department of Mathematics, Vrije Universiteit Amsterdam}
}

\maketitle

\begin{abstract}
It is natural to expect that there are only three possible types of scaling
limits for the collection of all percolation interfaces in the plane: (1) a
trivial one, consisting of no curves at all, (2) a critical one, in which all
points of the plane are surrounded by arbitrarily large loops and every deterministic
point is almost surely surrounded by a countably infinite family of nested loops
with radii going to zero, and (3) an intermediate one, in which every deterministic
point of the plane is almost surely surrounded by a largest loop and by a countably
infinite family of nested loops with radii going to zero. We show how one can
prove this using elementary arguments, with the help of known scaling relations for
percolation.

The trivial limit corresponds to subcritical and supercritical percolation,
as well as to the case when the density $p$ approaches the critical probability,
$p_c$, sufficiently slowly as the lattice spacing is sent to zero. The second
type corresponds to critical percolation and to a faster approach of $p$ to $p_c$.
The third, or near-critical, type of limit corresponds to an intermediate speed
of approach of $p$ to $p_c$. The fact that in the near-critical case a deterministic
point is a.s.\ surrounded by a largest loop demonstrates the persistence of a
macroscopic correlation length in the scaling limit and the absence
of scale invariance.
\end{abstract}



\section{Introduction and main results}

In Bernoulli site (respectively, bond) percolation, the sites
(resp., bonds) of a regular lattice with lattice spacing $\delta$
are colored white with probability $p$ and black otherwise,
independently of each other. One is then interested in the
connectivity properties of the monochromatic subgraphs of the
lattice, called \emph{clusters} (see, e.g., \cite{br,gr,kesten-book}).

The rigorous geometric analysis of the continuum scaling limit
($\delta \to 0$) of two-di\-men\-sional critical site percolation on
the triangular lattice has made tremendous progress in recent years.
In particular, the work of Schramm~\cite{schramm} and
Smirnov~\cite{smirnov} has allowed to identify the scaling limit of
critical interfaces (i.e., boundaries between black and white
clusters) in terms of the Schramm-Loewner Evolution (SLE) (see
also~\cite{cn3,smirnov06}). Based on that, Camia and Newman have
constructed~\cite{cn1} a process of continuum nonsimple loops in the
plane, and proved~\cite{cn2} that it coincides with the scaling
limit of the collection of all percolation interfaces (the
\emph{full} scaling limit). The use of SLE technology and
computations, combined with Kesten's scaling
relations~\cite{kesten87}, has also led to the derivation of
important properties of percolation such as the values of some
critical exponents~\cite{lsw5,sw}.

In later work~\cite{cfn1,cfn2}, based on heuristic arguments,
Camia, Fontes and Newman have proposed an approach for obtaining
a one-parameter family of \emph{near-critical} scaling limits with
density of white sites (or bonds) given by
\begin{equation} \label{eq-near-critical}
p=p_c+\lambda\delta^{\alpha},
\end{equation}
where $p_c$ is the critical density, $\delta$ is the lattice spacing,
$\lambda \in (-\infty,\infty)$, and $\alpha$ is set equal to $3/4$ to
get nontrivial $\lambda$-dependence in the limit $\delta \to 0$ (see
below and~\cite{aizenman,ab,bcks}). The approach proposed in~\cite{cfn1,cfn2}
is based on the critical full scaling limit and the ``Poissonian marking"
of some special (``macroscopically pivotal") points, and it leads to a
conceptual framework that can in principle describe not only the scaling
limit of near-critical percolation but also of related two-dimensional
models such as dynamical percolation, the minimal spanning tree and
invasion percolation (see~\cite{cfn2}).

In this note, we consider the collection of all percolation interfaces and show
how one can use known scaling relations for percolation to prove that, besides
the trivial scaling limit corresponding to subcritical and supercritical percolation,
there are only two other alternatives, that we call \emph{critical} and
\emph{near-critical} scaling limits, and for which we give a geometric characterization.


We postpone precise definitions till Sections \ref{sec-notation} and \ref{sec-scal-lim}
(including those of the space of interfaces and the topology of weak convergence),
but in order to avoid delaying the statement of the main result, we present it here.
We denote by $P_{\delta,p}$ the probability measure corresponding
to Bernoulli site percolation on the triangular lattice with lattice
spacing $\delta$ and parameter $p$. It is well known~\cite{kesten-book} that
percolation on the triangular lattice has a phase transition at $p=1/2$.
Let $H_{\delta}^w(n)$ denote the event that there is a white horizontal
crossing in a ``box'' of Euclidean side length $n\delta$ on the lattice
with lattice spacing $\delta$ (see Figure~\ref{fig-tribox} and the next
section for precise definitions). Due to the black/white symmetry of the
model, without loss of generality, we can restrict our attention to the
case $p \geq 1/2$. For $\epsilon \in (0,1/2)$, we define
\begin{eqnarray*}
p^+_{\epsilon}(n) := \inf\{p:P_{\delta,p}(H_{\delta}^w(n))>1/2+\epsilon\}.
\end{eqnarray*}
(Note that $p^+_{\epsilon}(n)$ is independent of $\delta$,
$p^+_{\epsilon}(n) \geq 1/2 \,\,\, \forall \epsilon \in (0,1/2)$ and
$p^+_{\epsilon_1}(n) \leq p^+_{\epsilon_2}(n)$ if
$\epsilon_1 \leq \epsilon_2$.)
Let $\mu_{\delta,p}$ denote the distribution of all percolation interfaces
for site percolation with parameter $p$ on the triangular lattice with lattice
spacing $\delta$.

\begin{theorem} \label{main-thm}
Suppose that $\mu$ is the weak limit of a sequence
$\{\mu_{\delta_j,p_j}\}_{j \in {\mathbb N}}$, with $\delta_j \to 0$ as
$j \to \infty$ and $p_j \geq 1/2$ for all $j$. Then one of the following
non-void scenarios holds.
\begin{enumerate}
\item[(1)] Trivial scaling limit: $\mu$-a.s.\ there are no loops of diameter larger than zero.
\item[(2)]  Critical scaling limit: $\mu$-a.s.\ any deterministic point in the plane is
surrounded by a countably infinite family of nested loops with radii
going to zero. Moreover, every point is surrounded by a countably infinite
family of nested loops with radii going to infinity.
\item[(3)] Near-critical scaling limit: $\mu$-a.s.\ any deterministic point in the plane is
surrounded by a largest loop and by a countably infinite family of
nested loops with radii going to zero.
\end{enumerate}
Moreover, the third scenario can be realized by taking $0<\epsilon_1<\epsilon_2<1/2$
and (an appropriate subsequence of) $\{p_j\}_{j \in {\mathbb N}}$ chosen so that
$p^+_{\epsilon_1}(1/\delta_j) \leq p_j \leq p^+_{\epsilon_2}(1/\delta_j)$
for every $j$.
\end{theorem}

The above geometric characterization of near-critical scaling
limits, case {\em(3)}, was conjectured in~\cite{cfn1}. It shows that
such limits are not scale invariant and differ qualitatively from
the critical scaling limit at large scales, since in the latter case
there is no largest loop around any point. At the same time, they
resemble the critical scaling limit at short scales because of the
presence, around any given point, of infinitely many nested loops
with radii going to zero. Depending on the context, this situation
is also described as {\em off-critical} or {\em massive} scaling
limit (where ``massive" refers to the persistence of a macroscopic
correlation length, which should give rise to what is known in the
physics literature as a ``massive field theory").

The three regimes in Theorem~\ref{main-thm} correspond to those in Proposition~4 of~\cite{nowe},
which contains, among other things, results analogous to some of ours in the context of
a single percolation interface and its scaling limit.
Perhaps the most interesting results of~\cite{nowe} and of this paper concern the
near-critical regime (regime {\em (3}) of Theorem~\ref{main-thm}). While~\cite{nowe}
deals with a single interface, proving that its scaling limit in the near-critical
regime is singular with respect to SLE$_6$ (the critical scaling limit), in this
paper we consider the full scaling limit and are concerned with the geometry of the
set of all interfaces, so that, in some sense, our result on the near-critical regime
complements that of~\cite{nowe}.

Our results imply that when $\alpha>3/4$ in~(\ref{eq-near-critical}),
the full scaling limit is trivial, when $\alpha<3/4$ it is critical, and there is a
non-empty regime where it is neither trivial nor critical. The following corollary is
an immediate consequence of Theorem~\ref{main-thm} and the power law~(\ref{powerlaw})
given at the end of the next section.

\begin{corollary} \label{cor-alpha}
Consider site percolation on the triangular lattice with lattice spacing $\delta$
and parameter $p=1/2+\lambda\delta^\alpha$. Then, for every $\lambda \in (-\infty,\infty)$,
\begin{itemize}
\item if $\alpha<3/4$, there is a unique scaling limit which is trivial in the sense
of Theorem~\ref{main-thm},
\item if $\alpha>3/4$, every subsequential scaling limit is critical in the sense
of Theorem~\ref{main-thm}.
\end{itemize}
\end{corollary}

It is natural to conjecture that the near-critical regime (case {\em (3)}
in Theorem~\ref{main-thm}) corresponds to the case $\alpha=3/4$, but at the
moment this is not known. In order to prove that, one would need to show
that when $\alpha=3/4$ the correlation length $L_{\epsilon}(p)$ defined in
Section~\ref{sec-notation} below remains bounded away from zero and infinity
as $\delta \to 0$ (see the proof of case {\em (3)} of Theorem~\ref{main-thm}).
This is believed to be the case, and in fact the correlation length is
expected to follow the power law $L_{\epsilon}(p) \asymp |p-1/2|^{-4/3}$
as $p \to 1/2$, where $\asymp$ means that the ratio between the two
quantities is bounded away from zero and infinity. However, only the
weaker power law $L_{\epsilon}(p) = |p-1/2|^{-4/3 + o(1)}$
has been proved~\cite{sw}.

In a remark at the end of Section~\ref{sec-proofs}, we explain how one can combine
\cite{cn2}  with (part of) the proof of Proposition 4 of~\cite{nowe} to obtain a
(much) stronger version of case {\em (2)} of Theorem~\ref{main-thm} (namely, that
the scaling limit in regime {\em (2)} coincides with the critical full scaling
limit~\cite{cn1,cn2}). In view of this result, in the second item of Corollary~\ref{cor-alpha},
one can identify the scaling limit with $\alpha>3/4$ with the unique critical scaling
limit~\cite{cn1,cn2}.

It is our understanding that significant progress has recently been made~\cite{gps}
(see also~\cite{garban}) in proving the approach of~\cite{cfn1,cfn2} to near-critical
scaling limits. A consequence would be that {\em all} subsequential limits discussed
in this paper are in fact limits.

To conclude this section, we point out that, although our results are stated for
site percolation on the triangular lattice, except for Corollary~\ref{cor-alpha}
which relies on the power law~(\ref{powerlaw}) and Remark 4.1 which relies on results
from~\cite{cn2,cn3,nowe}, they also apply to bond percolation and to other regular
lattices like the square lattice (after replacing $1/2$ with $p_c$ when necessary).
Indeed, the main tools in our proofs originated in Kesten's work~\cite{kesten87}
on the square lattice and can be used for both site and bond percolation models on
a large class of lattices (see~\cite{kesten-book}). For a discussion of the range
of applicability of Kesten's and related results, and consequently of the results
of the present paper, the reader is referred to Section~8.1 of~\cite{nolin}.

\section{Notation and some background} \label{sec-notation}

\begin{figure}[!ht]
\begin{center}
\includegraphics[width=5cm]{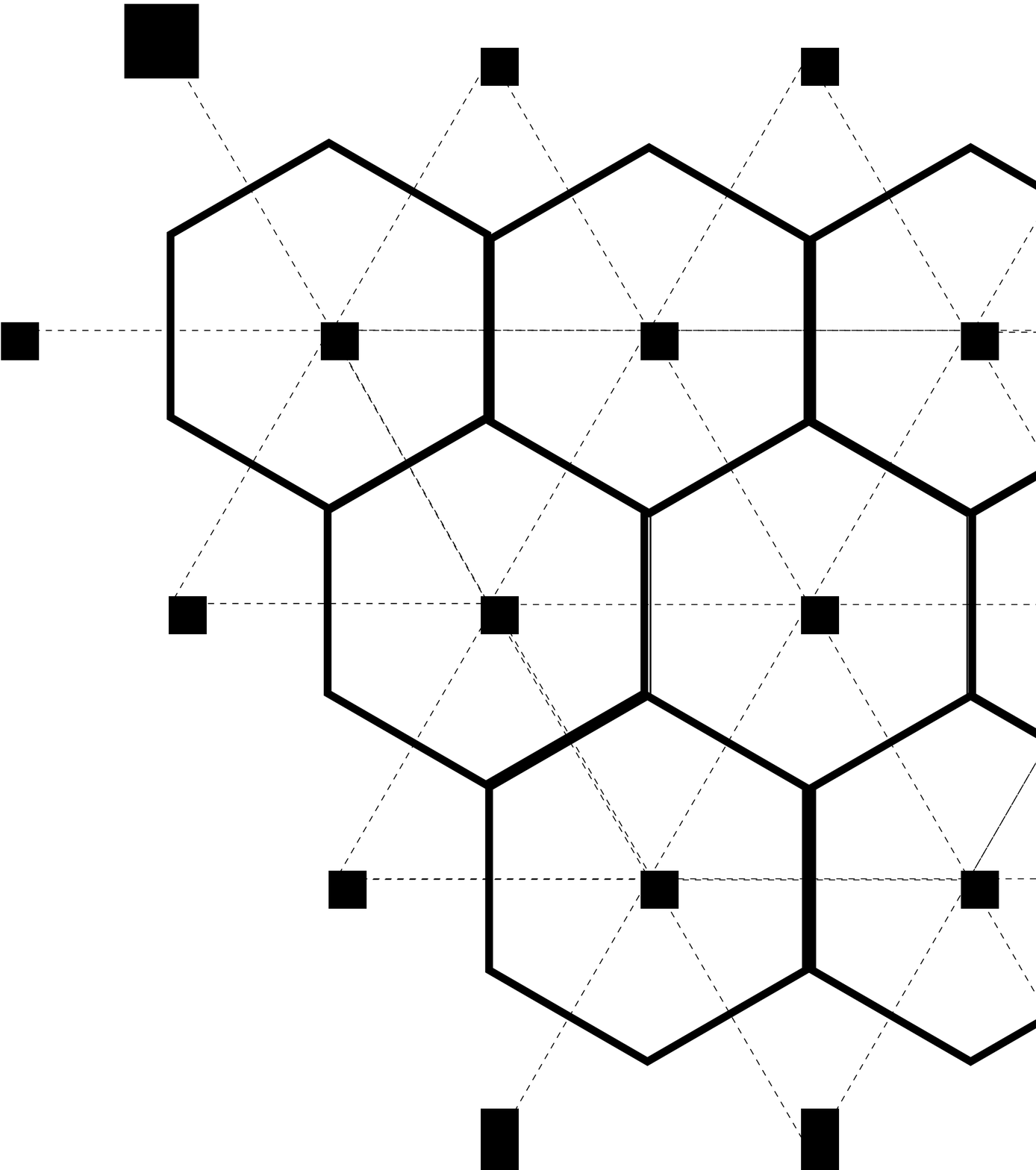}
\caption{Embedding of the triangular and hexagonal lattices in ${\mathbb R}^2$.}
\label{fig-tri-hex}
\end{center}
\end{figure}

Consider the hexagonal lattice $\HH_{\delta}$ with lattice
spacing $\delta>0$, and its dual, the triangular lattice
$\TT_{\delta}$, embedded in ${\IR}^2$ as in
Figure~\ref{fig-tri-hex}. A site of the triangular lattice is
identified with the face of the hexagonal lattice that contains it.

Throughout this paper, we are interested in Bernoulli site percolation
on $\TT_{\delta}$, defined as follows.
Each site of $\TT_{\delta}$ is independently declared white, and the
corresponding hexagon colored white, with probability $p$. Sites that
are not white are declared black, and the corresponding hexagons are
colored black. We denote by $P_{\delta,p}$ the probability measure
corresponding to site percolation on $\TT_{\delta}$ with parameter $p$.
It is well known~\cite{kesten-book} that percolation on the triangular
lattice has a phase transition at $p=1/2$.

A {\em path} of length $n$ in $\TT_{\delta}$ is a sequence of $n$ distinct
sites $(x_1,x_2,\ldots,x_n)$ of $\TT_{\delta}$ and the edges of $\TT_{\delta}$
between them such that $x_k$ and $x_{k+1}$ are adjacent in $\TT_{\delta}$ for
all $k=1,\ldots,n-1$. A {\em circuit} of length $n$ is a path whose
first and last sites are adjacent. We define the diameter of a set
$U \subset \mathbb{R}^2$ as
\[
\mbox{diam}(U):=\sup\{|x-y|:x,y \in U \},
\]
where $|\cdot|$ denotes Euclidean distance. We call a path or a circuit
white (resp., black) if all its sites are white (resp., black).

The edges between neighboring hexagons with different colors form
{\em interfaces}. A concatenation of such edges will be called a
{\em boundary path} or a {\em boundary loop} if it forms a closed
curve. Note that boundary curves and loops are always simple (i.e.,
no self-touching occurs) for $\delta>0$. However, this will
not necessarily be the case in the scaling limit $\delta \to 0$.

\begin{figure}[!ht]
\begin{center}
\includegraphics[width=5cm]{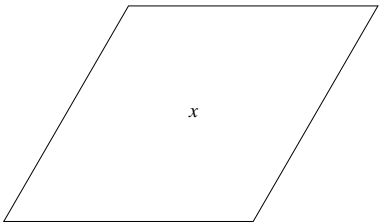}
\caption{The box $B(x;r)$.}
\label{fig-tribox}
\end{center}
\end{figure}

For $n_1,n_2>0$, $[0,n_1] \times [0,n_2]$ will denote the closed
parallelogram with Euclidean side-lengths $n_1$ and $n_2$ and sides
which are parallel to two of the axes of the triangular lattice
as in Figure~\ref{fig-tribox}. In particular, when $n_1=n_2$,
we call such a parallelogram a {\em box}. $B(x;r)$ will denote the
box centered at $x$, obtained by translating $[0,r] \times [0,r]$
(see Figure~\ref{fig-tribox}).
For $0 < r < R$, we define the {\em annulus} $A(x;r,R)$ as
\[
A(x;r,R):=B(x;R)\setminus B^{\circ}(x;r),
\]
where $B^{\circ}(x;r)$ denotes the interior of $B(x;r)$.
When $x$ is the origin, we will write $B(r)$ and $A(r,R)$, respectively.
Note that boxes and annuli are defined in terms of the Euclidean metric
and not relative to the lattice spacing.

The notion of {\em correlation length} will be very important.
Various equivalent definitions are possible; we choose the one,
introduced in~\cite{ccfs} and also used in~\cite{kesten87}, that is
most suitable for our purposes. Let $n$ be an integer and
$H_{\delta}^w(n)$ be the event that a percolation configuration on
$\TT_{\delta}$ contains a white path inside $B(n \delta)$
intersecting both its ``left side" and its ``right side." For each
$\epsilon \in (0,1/2)$, the correlation length $L_{\epsilon}(p)$ is
defined as follows:
\begin{eqnarray*}
L_{\epsilon}(p):= \min\{n:P_{\delta,p}(H_{\delta}^w(n))>1/2+\epsilon\} \mbox{ when
} p > 1/2,\\
L_{\epsilon}(p):= \min\{n:P_{\delta,p}(H_{\delta}^w(n))<1/2-\epsilon\} \mbox{ when
} p  < 1/2.
\end{eqnarray*}
We also define $L_{\epsilon}(1/2)=\infty$ for all $\epsilon \in (0,1/2)$.
Note that in the definition above, the correlation length is measured in
lattice spacings (rather than in the Euclidean metric), and is therefore
independent of $\delta$. Below we will frequently make use of the \emph{scaled}
correlation length $\delta L_{\epsilon}(p)$, which can be seen as the
``macroscopic'' correlation length.

An important fact about the correlation length is that the $\epsilon$ in
the definition is unimportant, due to the following result~\cite{nolin}
(a weaker version is proved in~\cite{kesten87}):
for any $\epsilon, \epsilon' \in (0,1/2)$ we have
\begin{equation} \nonumber
L_{\epsilon}(p) \asymp L_{\epsilon'}(p),
\end{equation}
where $f \asymp g$ means that the ratio between the functions $f$ and $g$ is
bounded away from 0 and $\infty$ as $p \to 1/2$. In view of this, we fix some
$\epsilon \in (0,1/2)$ and work with this choice of $\epsilon$ throughout the
rest of the paper without loss of generality.
We will also need the following five results. The first is a consequence
of Theorem~26 of~\cite{nolin} (see also Theorem~1 of~\cite{kesten87} for a
similar result).
\begin{lemma} \label{circuitjes}
Consider percolation on $\TT_{\delta}$ with parameter $p$ and let $C^w(r,R)$
(resp., $C^b(r,R)$) be the event that the annulus $A(r,R)$ is crossed (from
the inner to the outer boundary) by a white (resp., black) path.
Then,
$$
P_{\delta,p}(C^w(r,R)) = P_{\delta,p}(C^b(r,R)) \asymp P_{\delta,1/2}(C^b(r,R))
= P_{\delta,1/2}(C^w(r,R))
$$
uniformly in $p$ and $0 < r \leq R \leq \delta L_{\epsilon}(p)$.
\end{lemma}
We interpret this result as follows: on a scale not larger than the correlation
length, percolation with parameter $p$ looks roughly like critical percolation.

The second result, stated below, is Remark~38 of~\cite{nolin}.
\begin{lemma}\label{lemma remark nolin}
Consider percolation on $\TT_\delta$ with parameter $p\geq 1/2$. Let $C_{H}([0,n]\times[0,kn])$
denote the event that the parallelogram $[0,n]\times[0,kn]$ contains a \emph{black} horizontal
crossing. For any $k\geq1$ there exist two constants $C_1 < \infty$ and $C_2>0$, both depending
on $k$ and $\epsilon$, such that
\[
 P_{\delta,p}(C_{H}([0,n]\times[0,kn])) \leq C_{1} \exp\left(-\frac{C_{2} n}{\delta L_{\epsilon}(p)}\right).
\]

\end{lemma}

The third result is as follows (see, e.g., \cite{bpsv,nolin} for more explanation
and references).
\begin{lemma}\label{lemma Van den Berg}
Consider percolation on $\TT_{\delta}$ with parameter $p \geq 1/2$, and let
\[
D_r=\{\exists \mbox{ black circuit } S \mbox{ surrounding the origin with }
\emph{diam}(S) \geq  r\}.
\]
Then, for each $\epsilon \in (0,1/2)$ there exist two constants
$C_3=C_3(\epsilon)<\infty$ and $C_4=C_4(\epsilon)>0$ such that
\begin{equation} \nonumber 
P_{\delta,p}(D_r) \leq C_3\exp\left(-\frac{C_4 r}{\delta L_{\epsilon}(p)}\right).
\end{equation}
\end{lemma}

The fourth result is an immediate consequence of
Lemma~\ref{lemma remark nolin} (see Figure~\ref{fig-nocircuits}
for an example of a similar argument).

\begin{lemma}\label{lemma nolin}
Consider percolation on $\TT_{\delta}$ with parameter $p \geq 1/2$, and let
\[
D'_r=\{\exists \mbox{ black path containing the origin and of diameter at least }  r\}.
\]
For each $\epsilon \in (0,1/2)$ there exist two constants
$C_5=C_5(\epsilon)<\infty$ and $C_6=C_6(\epsilon)>0$ such that
\begin{equation} \nonumber 
P_{\delta,p}(D'_r) \leq C_5\exp\left(-\frac{C_6 r}{\delta L_{\epsilon}(p)}\right).
\end{equation}
\end{lemma}

The last result is the celebrated power law for the correlation length~\cite{sw}
(see also~\cite{nolin}): as $p \to 1/2$,
\begin{equation}
\label{powerlaw}
L_{\epsilon}(p) = |p-1/2|^{-4/3 + o(1)} .
\end{equation}

\section{The scaling limit} \label{sec-scal-lim}

We turn our attention to the main object of study in this article --
the scaling limit of the collection of all boundary loops. We will
follow the approach of~\cite{cn2}, using the topology introduced
in~\cite{ab}.


When taking the scaling limit as the lattice spacing
$\delta \to 0$ one can focus on fixed finite regions,
$\Lambda \subset {\mathbb R}^2$, or consider the whole
${\mathbb R}^2$ at once.
The second option avoids dealing with boundary conditions,
but requires an appropriate choice of metric.
A convenient way of dealing with the whole ${\mathbb R}^2$
is to replace the Euclidean metric with a distance function
$\Delta(\cdot,\cdot)$ defined on
${\mathbb R}^2 \times {\mathbb R}^2$ by
\begin{equation} \nonumber
\Delta(u,v) := \inf_{\varphi} \int (1 + | {\varphi(s)} |^2)^{-1} \, ds,
\end{equation}
where the infimum is over all smooth curves $\varphi(s)$
joining $u$ with $v$, parametrized by arclength $s$, and
where $|\cdot|$ denotes the Euclidean norm.
This metric is equivalent to the Euclidean metric in bounded
regions, but it has the advantage of making ${\mathbb R}^2$
precompact.
Adding a single point at infinity yields the compact space
$\dot{\mathbb R}^2$ which is isometric, via stereographic
projection, to the two-dimensional sphere.

In dealing with the scaling limit we use the approach of
Aizenman-Burchard~\cite{ab}.
We regard curves as equivalence classes of continuous
functions from the unit interval to $\dot{\mathbb R}^2$,
modulo monotonic reparametrizations. Below, $\gamma$ will
represent a particular curve and $\gamma(t)$ a parametrization
of $\gamma$. Denote by $\cal S$ the complete separable metric
space of curves in $\dot{\mathbb R}^2$ with the distance
\begin{equation} \label{Distance}
\text{D} (\gamma_1,\gamma_2) := \inf
\sup_{t \in [0,1]} \Delta(\gamma_1(t),\gamma_2(t)),
\end{equation}
where the infimum is over all choices of parametrizations
of $\gamma_1$ and $\gamma_2$ from the interval $[0,1]$.
A set of curves (more precisely, a closed subset of $\cal S$)
will be denoted by ${\cal F}$. The distance between two closed
sets of curves is defined by the induced Hausdorff metric as follows:
\begin{equation} \label{hausdorff-D}
\text{Dist}({\cal F},{\cal F}') \leq \varepsilon
\Leftrightarrow (\forall \, \gamma \in {\cal F}, \, \exists \,
\gamma' \in {\cal F}' \text{ with }
\text{D} (\gamma,\gamma') \leq \varepsilon
\text{ and vice versa}).
\end{equation}
The space $\Omega$ of closed sets of $\cal S$
(i.e., collections of curves in $\dot{\mathbb R}^2$)
with the metric~(\ref{hausdorff-D}) is also a complete
separable metric space. We denote by ${\cal B}$ its
Borel $\sigma$-algebra.

When we talk about convergence in distribution of random curves,
we always mean with respect to the uniform metric~(\ref{Distance}),
while when we deal with closed collections of curves, we always
refer to the metric~(\ref{hausdorff-D}). In this paper, the space
$\Omega$ of closed sets of $\cal S$ is used for collections of
boundary loops and their scaling limits.


Aizenman and Burchard \cite{ab} formulate a hypothesis that implies,
for every sequence $\delta_j \downarrow 0$, the existence of a
scaling limit along some subsequence $\{\delta_{j_i}\}$. The
hypothesis in \cite{ab} is formulated in terms of crossings of {\em
spherical} annuli, but one can work with the annuli defined in
Section~\ref{sec-notation} just as well. In order to state it, we
need one more piece of notation. For $\delta >0$, we denote by
$\mu_{\delta}$ any probability measure supported on collections of
curves that are polygonal paths on the edges of the hexagonal
lattice $\HH_{\delta}$.

In our context, the hypothesis is as follows.
\begin{hypothesis}\label{Hypothesis 1}
For all $k < \infty$ and for all annuli $A(x;r,R)$ with $\delta
\leq r \leq R\leq 1$, the following bound holds uniformly in
$\delta$:

\begin{equation} \nonumber
\mu_{\delta}\left(A(x;r,R) \mbox{ is crossed by } k \mbox{ disjoint curves}\right) \leq K_k
\left(\frac{r}{R}\right)^{\phi(k)}
\end{equation}
for some $K_k < \infty$ and $\phi(k) \rightarrow \infty$ as $k \rightarrow \infty$.
\end{hypothesis}

The next theorem follows from a more general result proved in~\cite{ab}.
\begin{theorem} [\cite{ab}] \label{Thm Aizenman}
Hypothesis \ref{Hypothesis 1} implies that for any sequence $\delta_j \downarrow 0$,
there exist a subsequence $\{\delta_{j_{i}}\}_{i \in {\mathbb N}}$ and a probability
measure $\mu$ on $\Omega$ such that $\mu_{\delta_{j_{i}}}$ converges weakly to $\mu$
as $i \rightarrow \infty$.
\end{theorem}

It was already remarked in the appendix of~\cite{ab} that the above hypothesis
can be verified for two-dimensional critical and near-critical percolation.
The same conclusion follows from results in~\cite{nolin},
and is obtained in Proposition~1 of~\cite{nowe}. We will need a slightly more
general result, stated and proved below for completeness.

\begin{lemma} \label{lemma-hypo}
Let $\{\mu_{\delta_j,p_j}\}_{j \in {\mathbb N}}$ be a sequence of measures on
boundary paths induced by percolation on $\TT_{\delta_j}$ with parameters $p_j$.
For {\em any} sequence $\delta_j \to 0$ and {\em any} choice of the collection
$\{p_j\}_{j \in {\mathbb N}}$, Hypothesis~\ref{Hypothesis 1} holds.
\end{lemma}

\medskip\noindent
{\bf Proof.} First of all, observe that the number of boundary paths
crossing an annulus is necessarily even and that, if there are $k$
disjoint boundary paths crossing the annulus, then the annulus must
also be crossed by $k/2$ disjoint black paths. For any
$\delta>0$ and $p \geq 1/2$,
we have
\begin{eqnarray*}
\lefteqn{P_{\delta,p}\left(A(x;r,R) \mbox{ is crossed by $k/2$ disjoint black paths} \right)} \\
& & \leq
P_{\delta,1/2}\left(A(x;r,R) \mbox{ is crossed by $k/2$ disjoint black paths} \right) \\
& & \leq P_{\delta,1/2}\left(A(x;r,R) \mbox{ is crossed by a black
path} \right)^{k/2},
\end{eqnarray*}
where we have used monotonicity and the BK inequality~\cite{vdbk}. Define
$l_{1}, l_2$ as the largest, respectively smallest integer such that $1/2^{l_{1}}\geq r$,
resp. $1/2^{l_{2}}\leq R$. Consider the annuli
\[
A_{1} = A(x;(1/2^{l_{1}},1/2^{l_{1}-1}), A_{2} = A(x;1/2^{l_{1}-1},1/2^{l_{1}-2}),
\ldots, A_{N} = A(x;1/2^{l_{2}+1},1/2^{l_{2}}),
\]
where $N$ denotes the maximal number of annuli of this type that can
be placed in $A(x;r,R)$. Note that $N$ is of order $\log(R/r)$ and
hence there exists a constant $C>0$, independent of $r$ and $R$,
such that $\lfloor N/2 \rfloor \geq C\log(R/r)$. Observe furthermore
that if $A(x;r,R)$ is crossed by a black path then none of the
annuli $A_{1},\ldots,A_{N}$ contains a white circuit surrounding
$x$. It follows from the RSW theorem (see, e.g.,
\cite{kesten-book,gr,nolin}) that the probability of the event that
the annulus $A_i$ contains a white circuit is uniformly (in $i$)
bounded from below by some $\gamma>0$, independent of
$\delta$. By definition of the annuli, $A_{2i+1}$ and $A_{2i^{'}+1}$
are disjoint for $i\neq i^{'}$. Putting everything together we
obtain
\begin{eqnarray*}
\lefteqn{P_{\delta,1/2}\left(A(x;r,R) \mbox{ is crossed by a black path} \right)^{k/2}} \\
& & \leq P_{\delta,1/2}\left(\bigcap_{i=0}^{\lfloor
N/2\rfloor-1}\{A_{2i+1} \text{ does not contain
a white circuit surrounding } x \} \right)^{k/2} \\
&& = \left[\prod_{i=0}^{\lfloor N/2\rfloor-1} P_{\delta,1/2}(A_{2i+1} \text{ does not contain
a white circuit surrounding } x)\right]^{k/2} \\
&& \leq \left[(1-\gamma)^{\lfloor N/2\rfloor}\right]^{k/2} \\
&& \leq \left[(1-\gamma)^{C \log(R/r)}\right]^{k/2} \\
&& = \left[(r/R)^{-C \log(1-\gamma)}\right]^{k/2}.
\end{eqnarray*}
Therefore, taking $K_k=1$ and $\phi(k)=-C \log(1-\gamma) k/2$, we
obtain a bound of the desired form since $-C\log(1-\gamma)>0$.

For any $p \leq 1/2$,
the same uniform bound follows from swapping white and black in the
above argument. Hence, we have the desired bound for any sequence
$\delta_j \to 0$ and any $\{p_j\}_{j \in {\mathbb N}}$ and the lemma
is proved.

\section{Proof of Theorem \ref{main-thm}} \label{sec-proofs}

We first show how assuming different behaviors for the correlation length
leads to the three scenarios described in the theorem. Later we will prove
that those three scenarios are non-void and are the only three possibilities. \\

\noindent{\em (1)} Suppose that for some $\epsilon \in (0,1/2)$,
$\delta_j L_{\epsilon}(p_j) \to 0$ as $j \to \infty$. Recall that
$\mu$ is the weak limit of a sequence $\{\mu_{\delta_{j},p_{j}}\}$
with $\delta_j \to 0$ as $j \to \infty$ and $p_{j}\geq 1/2$. Note
that it is actually the case that $p_{j}> 1/2$ for all but finitely
many $j$ since $\delta_{j} L_{\epsilon}(p_{j}) \to 0$.

\begin{figure}[!ht]
\begin{center}
\includegraphics[width=7cm]{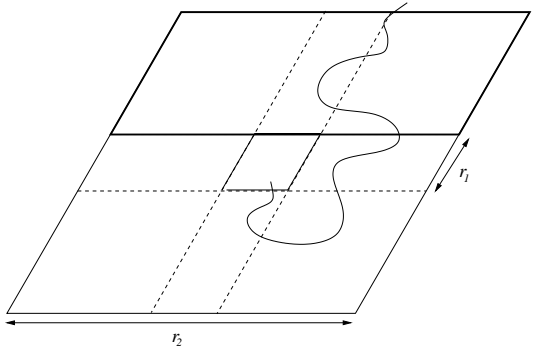}
\caption{The annulus $A(x;r_1,r_2)$ contains four rectangles of side lengths $(r_2-r_1)/2$ and $r_2$.
If a curve $\gamma$ crosses the annulus $A(x;r_1,r_2)$ then one (indicated with the heavy lines) of
the four rectangles contains a crossing in the ``easy'' direction. }
\label{fig-nocircuits}
\end{center}
\end{figure}

The existence of a boundary loop with positive diameter would imply that
there exist $x \in \IQ^2$ and $0< r_1 < r_2 \in \IQ$ such that the
loop intersects both $B^o(x,r_1)$ and ${\mathbb R}^2 \setminus B(x,r_2)$,
so that the annulus $A(x;r_1,r_2)$  is crossed by a boundary path. One of
the four (overlapping) parallelograms with side-lengths $(r_2-r_1)/2$ and
$r_2$ depicted in Figure \ref{fig-nocircuits} is then necessarily crossed
at least once in the ``easy'' direction by a boundary path (see Figure \ref{fig-nocircuits}).
Let $E$ denote such a crossing event. More precisely, crossings that realize
$E$ start and end outside the parallelogram and do not intersect the short
sides of the parallelogram. This makes $E$ open in our topology.
Note also that the occurrence of $E$ implies, for $\delta_j>0$, that the
parallelogram contains a black crossing in the easy direction. Thus, the
portmanteau theorem and Lemma~\ref{lemma remark nolin} yield
\begin{eqnarray}
\mu(A(x;r_1,r_2) \text{ is crossed by a boundary path}) & \leq& 4 \mu(E) \nonumber \\
& \leq& 4 \liminf_{j \to \infty} \mu_{\delta_{j},p_{j}}(E) \nonumber \\
& \leq& \liminf_{j \to \infty} C_{1} \exp \left(-\frac{ C_{2}(r_{2}-r_{1})}{2\delta_{j} L_{\epsilon}(p_{j})}\right) = 0. \nonumber
\end{eqnarray}
We can then conclude that
\begin{eqnarray*}
\lefteqn{\mu(\text{there exists a boundary loop with positive diameter})}\\
&&  \leq \bigcup_{x \in \IQ^{2}; r_{1},r_{2} \in \IQ^{+}} \mu(A(x;r_1,r_2) \text{ is
crossed by a boundary path}) =0.
\end{eqnarray*}

\noindent {\em (2)} Suppose that for some $\epsilon \in (0,1/2)$,
$\delta_j L_{\epsilon}(p_j) \to \infty$ as $j \to \infty$. Then
$1 \leq \delta_j L_{\epsilon}(p_j)$ for each $j$ sufficiently large.
To show the a.s.\ existence of an infinite sequence of boundary
loops with radii going to zero around the origin (or any other deterministic
point), we proceed as follows. Consider the sequence of annuli
$A_1=A(\frac{1}{2},1), A_2=A(\frac{1}{4},\frac{1}{2}), \ldots, A_k=A((\frac{1}{2})^k,(\frac{1}{2})^{k-1}), \ldots$ ,
and denote by $F_k$ the event that there is (at least) one boundary loop
surrounding the origin in the annulus $A_{2k+1}$ with $k\geq0$. Since we
defined annuli to be closed sets, the event $F_k$ is closed in our topology.
Note that in order to guarantee the presence of a boundary loop inside the
annulus $A_k$, it suffices to have, for example, a white circuit in
$A(\frac{5}{3} (\frac{1}{2})^k,(\frac{1}{2})^{k-1})$ and a black circuit in
$A((\frac{1}{2})^k,\frac{1}{3} (\frac{1}{2})^{k-2})$ (note that those two annuli
are disjoint). Since $\delta_j L_{\epsilon}(p_j) \geq 1$ for each large $j$, it
follows from Lemma~\ref{circuitjes} and the RSW theorem that the probability to
find a white circuit in $A(\frac{5}{3} (\frac{1}{2})^k,(\frac{1}{2})^{k-1})$ and
a black circuit in $A((\frac{1}{2})^k,\frac{1}{3} (\frac{1}{2})^{k-2})$ is bounded
away from $0$ as $j \to \infty$, uniformly in $k$. Therefore there exists $\varepsilon_0>0$
such that
\[\
\mu(F_k) \geq \limsup_{j \to \infty}\mu_{\delta_j,p_j}(F_k) \geq \varepsilon_0, \quad \text{for every } k,
\]
where the first inequality follows from the portmanteau theorem. Note also that the events $F_k$ and $F_{k'}$ are independent for $k' \neq k$. Moreover,
$\sum_{k=0}^{\infty} \mu(F_k) = \infty$ and thus by the Borel-Cantelli lemma there are
infinitely many boundary circuits surrounding the origin with diameter going to zero, $\mu$-a.s.

We argue in a similar way as above to show that every point is surrounded by a countably infinite
family of nested loops with radii going to infinity.
Let $B_k$ denote the annulus $A(2^{k},2^{k+1})$ and write $F'_k$, with $k\geq 0$, for the event that
the annulus $B_{2k+1}$ contains at least one boundary loop surrounding $B(1)$. For each $k$
it holds that $\delta_j L_{\epsilon}(p_{j})\geq 2^{k+1}$ for $j$ sufficiently large. Hence, it follows
from Lemma~\ref{circuitjes} and the RSW theorem that the probability to find a boundary loop in
$B_{2k+1}$ is uniformly (in $k$) bounded away from 0 as $j \to \infty$. Again, the event $F'_k$
is closed in our topology, thus
$\mu(F'_k) \geq \limsup_{j\to \infty}\mu_{\delta_{j},p_{j}}(F'_k) \geq \varepsilon_1$, for some
$\varepsilon_1>0$ independent of $k$. Hence, the Borel-Cantelli lemma implies that there
are infinitely many boundary loops surrounding $B(1)$ with diameter going to infinity
$\mu$-a.s. By translation invariance, the same is true for every $B(x;1)$
with $x \in {\mathbb Q}$ and therefore for every point of the plane. \\

\noindent{\em (3)} Suppose that for some $\epsilon \in (0,1/2)$,
$\delta_j L_{\epsilon}(p_j)$ stays bounded away from both $0$ and $\infty$
as $j \to \infty$. That is, there exist $\beta>0$ and $K<\infty$ such that
$\beta \leq \delta_j L_{\epsilon}(p_j) \leq K$ for each $j$ sufficiently large.
The first part of the proof in case {\em (2)} carries over directly to the present
case, with 1 replaced by $\beta$ in the lower bound for the macroscopic correlation
length and the annuli $A(1/2^k,1/2^{k-1})$ replaced by $A(\beta/2^k,\beta/2^{k-1})$.
Thus, $\mu$-a.s. there exist infinitely many boundary loops surrounding the origin,
with diameter going to zero.

Our next goal is to prove the a.s.\ existence of a largest boundary loop
surrounding the origin. Let $G_L$ denote the event that there exists a largest
boundary loop $\gamma$ surrounding or containing the origin and that this loop
has $\text{diam}(\gamma) \leq L$.
Then $G := \bigcup_{L=1}^{\infty} G_L$ is the event that there exists
a largest loop surrounding or containing the origin. Note that if all
black circuits around the origin have diameter smaller than $L$ and
there is no black path containing the origin of diameter larger than
$L-2\delta_{j}$, then $G_L$ occurs. Therefore
$\mu_{\delta_j,p_j}(G_L) \geq 1 - [P_{\delta_j,p_j}(D_L) + P_{\delta_j,p_j}(D'_{L-2\delta_{j}})]$,
where $D_L$ is the event that the origin is surrounded by a black circuit of diameter
at least $L$ and $D'_{L-2\delta_{j}}$ is the event that there is a black path
containing the origin of diameter at least $L-2\delta_{j}$. Using Lemmas \ref{lemma Van den Berg}
and \ref{lemma nolin} and the fact that the event $G_L$ is closed in our topology, we can write
\begin{eqnarray*}
\mu(G_L) &\geq& \limsup_{j \rightarrow \infty}
\mu_{\delta_{j},p_j}(G_L)  \nonumber \\
&\geq& 1-\liminf_{j\rightarrow \infty} [P_{\delta_j,p_j}(D_L) + P_{\delta_j,p_j}(D'_{L-2\delta_{j}})] \\
&\geq& 1-\liminf_{j \to \infty} C' \exp\left(-\frac{C'' L}{\delta_j L_{\epsilon}(p_j)}\right) \\
&\geq& 1- C' \exp\left(-\frac{C'' L}{K}\right). \label{largest-circuit}
\end{eqnarray*}
Since the events are nested (i.e., $G_{L_1} \subset G_{L_2}$ for $L_1<L_2$),
\[
\mu(G) = \lim_{L \rightarrow \infty} \mu(G_L)
\geq \lim_{L \rightarrow \infty} \left[1- C^{'} \exp\left(-\frac{C^{''} L}{K}\right)\right] = 1.
\]

Since boundary loops cannot cross each other and, by the previous part
of the proof, the origin is surrounded with probability one by a sequence
of infinitely many boundary loops with diameter going to zero, the largest
boundary loop does not touch the origin. Hence, the event $G$ coincides with
the existence of a largest boundary circuit surrounding the origin and we are done.

\bigskip

To continue the proof, note that for each $\epsilon \in (0,1/2)$, as
$j \to \infty$,
\begin{itemize}
\item either $\delta_j L_{\epsilon}(p_j) \to 0$,
\item or $\delta_j L_{\epsilon}(p_j) \to \infty$ ,
\item or $\delta_j L_{\epsilon}(p_j)$ is bounded away from both $0$
and $\infty$.
\end{itemize}
This is clearly so because we are assuming that
$\{\mu_{\delta_j,p_j}\}_{j \in {\mathbb N}}$ has a limit $\mu$, and we
have proved that the three cases above give rise to three incompatible
scenarios for $\mu$. Indeed, if we are not in one of the three cases
above, then there must be two different subsequences of
$\{(\delta_j,p_j)\}_{j \in {\mathbb N}}$ falling in two different cases,
which contradicts the existence of a limit $\mu$. We can then conclude
that there are no other possible scenarios for $\mu$ besides the three
described in the theorem.


To conclude the proof, we need to show that all three scenarios are non-void.
For the first two, this is obvious. To prove that the third scenario is also
non-void, take $0<\epsilon_1<\epsilon_2<1/2$ and consider any sequence
$\{(\delta_j,p_j)\}_{j \in {\mathbb N}}$ such that $\delta_j \to 0$ and
$p^+_{\epsilon_1}(1/\delta_j) \leq p_j \leq p^+_{\epsilon_2}(1/\delta_j)$.
This implies that
$L_{\epsilon_1}(p_j) \leq 1/\delta_j \leq L_{\epsilon_2}(p_j)$ for each $j$.
We can assume without loss of generality that the sequence
$\{\mu_{\delta_j,p_j}\}_{j \in {\mathbb N}}$ has a limit $\mu$. (If that
is not the case, by Theorem~\ref{Thm Aizenman} and Lemma~\ref{lemma-hypo}
we can extract a subsequence $\{\mu_{\delta_{j_k},p_{j_k}}\}_{k \in {\mathbb N}}$
that does have a limit, and rename it $\{\mu_{\delta_j,p_j}\}_{j \in {\mathbb N}}$.)
Since $\delta_j L_{\epsilon_1}(p_j)$ remains bounded as $j \to \infty$,
$\delta_j L_{\epsilon}(p_j)$ must remain bounded as $j \to \infty$
for every other $\epsilon \in (0,1/2)$ because
$L_{\epsilon}(p) \asymp L_{\epsilon_1}(p)$. Analogously, since
$\delta_j L_{\epsilon_2}(p_j)$ is bounded away from $0$ as $j \to \infty$,
$\delta_j L_{\epsilon}(p_j)$ must be bounded away from $0$ as $j \to \infty$
for every $\epsilon \in (0,1/2)$. Therefore, for each $\epsilon \in (0,1/2)$,
$\delta_j L_{\epsilon}(p_j)$ remains bounded away from both $0$ and
$\infty$ as $j \to \infty$, showing that $\mu$ falls in the third scenario. \\

\noindent\textbf{Remark 4.1.}
A significantly stronger version of case {\em (2)} can be proven, namely that
$\mu$ in this case coincides with the full scaling limit of critical percolation.
In order to prove this one can use the same strategy as in \cite{cn2}, combined
with the proof of Proposition 4 of \cite{nowe}.
Below we briefly sketch how one can obtain the result by modifying the arguments
of~\cite{cn2}. We stress that this is not meant to be a self-contained proof, and
some familiarity with~\cite{cn2} is needed in order to follow the arguments outlined
below.

First of all, wherever statement (S) (concerning the convergence of the critical
exploration path to SLE$_6$ --- see p.~18 of~\cite{cn2}) is invoked in~\cite{cn2},
one needs to use the proof of Proposition~4 of~\cite{nowe}. We remark that the
statement of Proposition 4 of~\cite{nowe} concerns triangular domains and is therefore
not sufficient for our purposes, but the proof applies in much greater generality.
In particular, the reader can check that it applies to the situations that arise
in~\cite{cn2}.

Uniform bounds on ``six-arm" events in the plane and ``three-arm" events near a
boundary (see Lemma~6.1 of~\cite{cn2} and its proof) are used repeatedly in~\cite{cn2}.
Thanks to Theorem~1 of~\cite{kesten87} and similar results described in~\cite{nolin}
(see, e.g., Section~3.2, and Theorem~27 and the discussion following it),
such uniform bounds are also available for percolation with parameter $p_j$
on ${\cal T}_{\delta_j}$ inside a disc of diameter $L$, provided that
$\delta_j L_{\epsilon}(p_j) \geq L$. Due to the assumption that
$\delta_j L_{\epsilon}(p_j) \to \infty$ as $j \to \infty$, such a condition
is satisfied for any $L<\infty$, for $j$ sufficiently large.

Some care is also needed in the proof of the second part of Theorem~5 of~\cite{cn2}
(see p.~19 for the statement; the proof begins on p.~27), and in particular of
Lemma~6.4 (see p.~27) and Lemma~6.6 (see p.~29) used in that proof. The proof
of the second part of Theorem~5 of~\cite{cn2} is given for critical percolation,
but the only features of critical percolation that are really used are uniform
bounds on certain crossing probabilities (involving crossings of a rectangle or
an annulus). Due to the above considerations, similar bounds can be used in the
present context (see again Section~3.2 and Theorem~27 of~\cite{nolin}).

\bigskip

\noindent {\bf Acknowledgements.} The first author thanks C.M. Newman for a useful discussion.

\end{document}